# A COMBINED LOGARITHMIC BOUND ON THE CHROMATIC INDEX OF MULTIGRAPHS

Michael Plantholt

Department of Mathematics, Illinois State University

Normal, IL 61790-4520  USA

Mikep@ilstu.edu ; (309) 438-7174

## Abstract

For a multigraph G, the integer round-up $\phi(G)$ of the fractional chromatic index $\chi_f(G)$ yields a good general lower bound for the chromatic index $\chi'(G)$. For an upper bound, Kahn showed that for any real $c > 0$ there exists a positive integer $N$ so that $\chi'(G) < \chi_f(G) + c\,\chi_f(G)$ whenever $\chi_f(G) > N$. We show that for any multigraph $G$ with order $n \geq 3$ and at least one edge, $\chi'(G) \leq \phi(G) + \log_{3/2}(\min\{(n+1)/3, \phi(G)\})$. Kahn's result follows as a corollary.

# A COMBINED LOGARITHMIC BOUND ON THE CHROMATIC INDEX OF MULTIGRAPHS

Michael Plantholt

Department of Mathematics, Illinois State University

Normal, IL 61790-4520  USA

## Abstract

For a multigraph $G$, the integer round-up $\phi(G)$ of the fractional chromatic index $\chi_f(G)$ yields a good general lower bound for the chromatic index $\chi'(G)$. For an upper bound, Kahn showed that for any real $c > 0$ there exists a positive integer $N$ so that $\chi'(G) < \chi_f(G) + c\,\chi_f(G)$ whenever $\chi_f(G) > N$. We show that for any multigraph $G$ with order $n \geq 3$ and at least one edge,
$\chi'(G) \leq \phi(G) + \log_{3/2}(\min\{(n+1)/3, \phi(G)\})$. Kahn's result follows as a corollary.

## 1. The Result

The *order* and *size* of a multigraph $G$ with vertex set $V(G)$ and edge set $E(G)$ are the cardinalities of $V(G)$ and $E(G)$, respectively, and denoted $n(G)$ and $e(G)$. Otherwise, we follow the terminology and notation of [4].

Recall that the chromatic index $\chi'(G)$ of a multigraph $G$ is the minimum number of matchings that are required to cover the edges of $G$. Clearly $\chi'(G)$ is at least as large as the maximum degree $\Delta(G)$. Moreover, if $H$ is any multigraph with odd order at least three, then $\chi'(H) \geq 2e(H)/(n(H) - 1)$, because any matching in $H$ can contain at most $(n(H) - 1)/2$ edges; we denote this lower bound on $\chi'(H)$ by $t(H)$. For any subset $S$ of the vertices of $G$ we let $<S>$ denote the subgraph of $G$ induced by the vertices in $S$. Then because $\chi'(G) \geq \chi'(H)$ for any subgraph $H$ of $G$, we have another lower bound for $\chi'(G)$ given by $\Gamma(G) = \max\{t(<S>)\}$, where the maximum is taken over all subsets $S$ of $V(G)$ for which $|S|$ is odd and at least 3.

Research supported in part by Summer Research Grants from
Illinois State University

Combining the lower bounds from the previous paragraph, we get a combined lower bound $\chi_f(G) = \max\{\Delta(G), \Gamma(G)\}$ for the chromatic index of G, and $\chi_f(G)$ is commonly referred to as the *fractional chromatic index* of G (a more complete discussion of this invariant is provided in [12]). When $n(G) < 3$, $\Gamma(G)$ is considered undefined and we let $\chi_f(G) = \Delta(G)$. Finally, because the chromatic index is always an integer, we get the standard combined lower bound $\chi'(G) \geq \phi(G)$, where $\phi(G)$ denotes the integer round-up of $\chi_f(G)$. Goldberg [2], [3] and Seymour [11] independently conjectured that this lower bound is quite tight, in the following sense (Goldberg's Conjecture was a bit stronger than the version stated here).

**Conjecture A**. For any multigraph G, $\chi'(G) \leq \max\{\Delta(G) + 1, \lceil \Gamma(G) \rceil\}$.

We often find it convenient to work with the following slightly weaker form of Conjecture A; this form of the Conjecture also appeared in [11].

**Conjecture B.** For any multigraph G, $\chi'(G) \leq 1 + \max\{\Delta(G), \lceil \Gamma(G) \rceil\}$.

Although Conjecture B is somewhat weaker than Conjecture A, it still reduces the possibilities for $\chi'(G)$ to two: either $\chi'(G) = \phi(G)$ or $\chi'(G) = \phi(G) + 1$.

The following three theorems give upper bounds for the chromatic index in terms of different invariants.

**Theorem A.** (Vizing [14]) If G is a simple graph, $\chi'(G) \leq \Delta(G) + 1$. More generally, for any multigraph G with maximum edge multiplicity m, $\chi'(G) \leq \Delta(G) + m$.

The next result is the most recent in a string of similar earlier results by Shannon [13], Andersen [1] and Goldberg [2], [3].

**Theorem B.** (Nishizeki & Kashiwagi [7]) Let G be a multigraph. If $\chi'(G) > (11\Delta(G) + 8)/10$, then $\chi'(G) = \phi(G) = \lceil \Gamma(G) \rceil$.

**Theorem C.** (Plantholt [9]) For any multigraph G of order n, $\chi'(G) \leq \phi(G) + \lceil n/8 \rceil - 1$.

In the three theorems above, we are given a bound for the amount by which the chromatic index can surpass a standard lower bound for $\chi'(G)$ in terms of three different invariants: the maximum edge multiplicity, the maximum degree, and the order of the multigraph. However, the amount by which the upper bound exceeds the lower bound is in each case linear in the chosen invariant, even though Conjecture B states that the difference between $\chi'(G)$ and $\phi(G)$ should be bounded by a constant (indeed, the constant 1). Using breakthrough methods, Jeff Kahn was able to provide evidence that improvements in this direction are possible by obtaining the following result:

**Theorem D.** (Kahn [5]). For any real $c > 0$ there exists a positive integer $N$ so that $\chi'(G) < \chi_f(G) + c\,\chi_f(G)$ whenever $\chi_f(G) > N$.

Using more basic methods, some similar type results were obtained in [10], but using the order $n$ as the key invariant:

**Theorem E.** For any multigraph $G$ with even order $n \geq 572$,
$\chi'(G) \leq \phi(G) + 1 + \sqrt{(n \ln n / 10)}$.

**Theorem F.** For any real $c > 0$, there exists a positive integer $N$ such that $\chi'(G) < \phi(G) + cn$ for any multigraph $G$ with order $n > N$.

We extend the methods from [10] to improve the results in two ways. First, the amount by which the upper bound surpasses $\phi(G)$ is now strictly logarithmic. Also, as the key parameter we now use the combined term $\min\{(n+1)/3, \phi(G)\}$, so that we get essentially two different upper bounds, choosing the minimum. Specifically, we show the following.

**Theorem 1.** For any multigraph $G$ with order $n \geq 3$ and at least one edge,
$\chi'(G) \leq \phi(G) + \log_{3/2}(\min\{n/3, \phi(G)\})$ if $n$ is even, and
$\chi'(G) \leq \phi(G) + \log_{3/2}(\min\{(n+1)/3, \phi(G)\})$ if $n$ is odd.

Kahn's result in Theorem D and the result in Theorem F then both follow as immediate corollaries.

## 2. The Approach

Let $G$ be a multigraph and let $S$ be a subset of $V(G)$. We let $<S>$ denote the subgraph of $G$ that is induced by the vertices in $S$, and let $\partial(S)$ denote the coboundary of $S$, that is, the cutset consisting of all edges that are incident with exactly one vertex of $S$. Recall that $\Delta(<S>)$ denotes the maximum degree of $<S>$. If we wish to stress that $<S>$ is being considered as an induced subgraph of $G$, we use the notation $\Delta(<S>; G)$ in order to clarify the host multigraph. We use the same convention for other invariants. For example, if $F$ is a 1-factor of $G$, then $\partial(S; G-F)$ gives the coboundary of the set of vertices $S$ within $G-F$.

Our approach is based on two elementary operations. One, removal of a matching, corresponds to assigning a color to the edges of the matching, and then eliminating them from further consideration. The second operation is the splitting operation, based on a shrinking procedure; these are defined below.

Let $G$ be a multigraph, and let $S$ be a non-empty proper subset of $V(G)$. The multigraph $G_S$ with vertex and edge sets described below is called the *multigraph obtained from $G$ by shrinking $S$*:

1. The vertex set of $G_S$ is given by $V(G_S) = V(G) - S \cup \{s\}$, where $s$ is the *vertex replacing $S$*.
2. Each vertex $u \neq s$ in $V(G_S)$ has exactly $|\partial(u,S)|$ edges in $G_S$ joining it to $s$, where $\partial(u,S)$ denotes the set of edges in $G$ that join $u$ with vertices of $S$. These edges are the only edges of $G_S$ that are incident with $s$, and so these comprise $\partial(\{s\})$ completely in $G_S$.
3. The edge set of $G_S$ is given by $E(G_S) = E(G-S) \cup \partial(\{s\})$.

For a subset $S$ of the vertex set of a multigraph $G$, we let $S^c$ denote the complementary set $V(G) - S$. We call the process of replacing $G$ by the pair of multigraphs $G_S$, $G_{S^c}$ a *splitting of $G$*. The proof of the main result will follow from a combination of 1-factor (or near 1-factor) removals and splittings. We depict the process using a binary tree, with root vertex $G$. When a vertex in the tree has only one child, the edge leading to that child represents the removal of a number of matchings; that number is given as the edge weight. A splitting is denoted by a vertex having two children, the resulting multigraphs. An example is given in Figure 1, which has starting multigraph G, and terminal multigraphs $Z_1$ through $Z_5$. In that example, there are $k_1$ matching removals before the first splitting. The sum of the edge weights along the path from the root to a terminal multigraph gives the total number of matchings removed during only the steps represented by that path; we call this number for a terminal multigraph Z its *matching removal number*, denoted mr(Z). For example, in Figure 1 we have mr($Z_3$) = $k_1$ + $k_2$ + $k_4$. Also, if $Z_i$ is a multigraph obtained at some point during this procedure, we call the value

mr($Z_i$) + $\phi$($Z_i$) - $\phi$(G)  the *cost* of going from G to $Z_i$. This cost we denote by cost(G → $Z_i$) ; it gives a measure of the number of "wasted" colors when matching removals did not reduce the $\phi$ value.

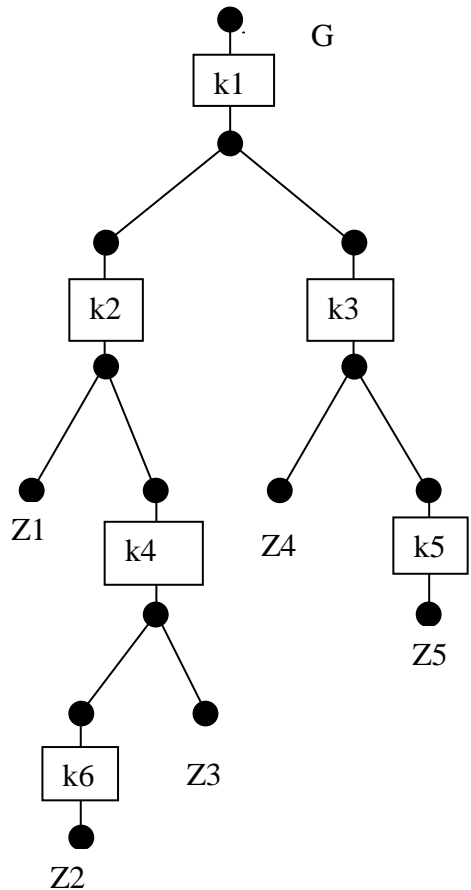

Figure 1. Diagram representing reduction process

In our main algorithm, our splittings will never increase the current ϕ-value. And by using splittings to control the multigraphs being considered appropriately, we are able to ensure that removal of a matching will in many cases (though of course not always) reduce the ϕ-value by 1. Since we can do this and significantly reduce the size (here measured by min { n/3 , ϕ }) of the terminal multigraphs, the result will follow. As a simple example, consider Figure 2 below. The multigraph G in (a) has n = Δ = ϕ = 6. We remove a 1-factor from G to obtain the multigraph in (b), which is 5-regular but with ϕ=6 still. Based on the overfullness in the triangles, we shrink on them to get the splitting shown in (c). The multigraphs in (c) still have ϕ=6, but since each has order n = 4, we've reduced min { n/3, ϕ }, multiplying it in each by a factor of 2/3; this has been achieved at the "cost" of one color (edge-matching removal) since a 1-factor has been removed without decreasing ϕ.

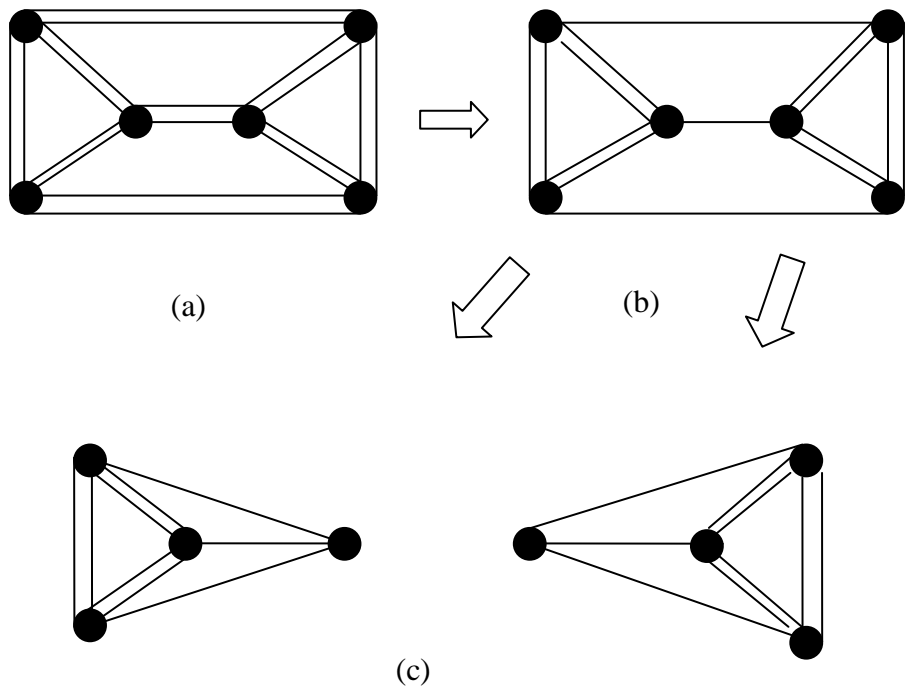

(a)        (b)

(c)

Figure 2. Reduction example.

The following well-known result uses the splitting procedure to bound the chromatic index.

**Lemma A.** ([11]) Let $G$ be a multigraph, and let $S$ be a non-empty proper subset of $V(G)$. Then $\chi'(G) \leq \max\{\chi'(G_S), \chi'(G_{S^c})\}$.

It is worth noting that given proper edge-colorings of $G_S$ and $G_{S^c}$ in $k = \max\{\chi'(G_S), \chi'(G_{S^c})\}$ colors, it is easy (see [11]) to obtain a corresponding proper k-coloring of the edges of $G$. Thus the bound we obtain by our algorithm could be applied to actually construct an edge-coloring. We also note that the fractional chromatic index version of Lemma A is also true, that is $\chi_f(G) \leq \max\{\chi_f(G_S), \chi_f(G_{S^c})\}$.

We will actually use the following extension of Lemma A, which follows immediately from Lemma A.

**Lemma A'.** Suppose that starting with multigraph $G$, we perform removals of matchings and $t$ splittings to reach terminal multigraphs $Z_1, Z_2, \ldots, Z_{t+1}$. Then
$\chi'(G) \leq \max\{mr(Z_1) + \chi'(Z_1), mr(Z_2) + \chi'(Z_2), \ldots, mr(Z_{t+1}) + \chi'(Z_{t+1})\}$.

By a *nontrivial multigraph* we mean one with more than one vertex, and a *nontrivial splitting* shrinks multigraphs that are both nontrivial. Recall that if $H$ is any multigraph with odd order at least three, we let $t(H) = 2e(H)/(n(H) - 1)$. For any integer $k$, a non-trivial odd order subgraph $H$ of $G$ is said to be *k-overfull* if $t(H) > k$, and *k-full* if $t(H) = k$; typically, we will be interested in fullness in the case where $k = \Delta(G)$. Recall also that $\Gamma(G) = \max\{t(<S>)\}$, where the maximum is taken over all nontrivial odd cardinality subsets $S$ of $V(G)$, the fractional chromatic index of $G$ is given by $\chi_f(G) = \max\{\Delta(G), \Gamma(G)\}$, and that the lower bound $\phi(G)$ for the chromatic index is given by $\phi(G) = \max\{\Delta(G), \lceil \Gamma(G) \rceil\}$. Finally, following terminology introduced in [11], for $r > 0$ we define an *r-graph* to be an r-regular multigraph $G$ for which $\Gamma(G) \leq r$ (in fact, if $G$ is an r-graph, then $G$ must have even order, and we must have $\Gamma(G) = r$ because for any vertex $v$, $t(G - v) = r$).

Because we can always add an isolated vertex to a multigraph, we need only prove Theorem 1 for even order multigraphs, that is

$\chi'(G) \leq \phi(G) + \log_{3/2}(\min\{n/3, \phi(G)\})$ if n is even.

The following result shows that we can also restrict our attention to r-graphs.

**Lemma B.** (Seymour [11]) Let G be any multigraph, and let $\phi(G) = r$. There exists a multigraph $G^+$ containing G such that
  i) $G^+$ is an r-graph, and
  ii) $|V(G^+)| = |V(G)|$ if $|V(G)|$ is even,
     $|V(G^+)| = |V(G)| + 1$ if $|V(G)|$ is odd.

Our proof is broken into three parts. In Part 1, we give Stage 1 of our algorithm. It begins by taking as input an r-graph on an even number n of vertices. By iteratively removing 1-factors and splitting, we construct larger and larger decomposition trees for G, without ever exceeding cost 1, until each branch results in a multigraph H satisfying one of the following halting conditions:

   Criterion 1A. $n(H) \leq 8$ or $\Delta(H) \leq 2$ ; $\text{cost}(G \to H) \leq 1$.
   Criterion 1B. $n(H) \leq (2/3) n(G)$, $\phi(H) \geq n(H)/2$, and $\text{cost}(G \to H) \leq 1$.
   Criterion 1C. (Termination for Stage 1 only; must proceed to Stage 2, which will involve removing near-1-factors) Min degree(H) = 0, $\Delta(H) \leq n(G)/3$.

The key result in Part 1 is Proposition 1, which shows that we can iterate the procedure until reaching the halting criteria without exacting cost greater than 1.

In Part 2, we present Stage 2 of the algorithm, which must be used when in Stage 1 we reach Criterion 1C. The Stage 2 halting criteria are:

   Criterion 2A. $n(H) \leq 8$ or $\Delta(H) \leq 2$ ; $\text{cost}(G \to H) \leq 1$.
   Criterion 2B. $n(H) \leq (2/3) n(G)$ , $\phi(H) \geq n(H)/2$, $\text{cost}(G \to H) \leq 1$.
   Criterion 2C. $n(H) \leq \phi(G)$ , $n(H) \leq n(G)/3$ , $\text{cost}(G \to H) \leq 2$ .
   Criterion 2D. $\phi(H) \leq 2 n(G)/9$, $\phi(H) \leq \phi(G)/2$, $\text{cost}(G \to H) \leq 1$.

Stage 2 is presented in algorithmic format, with arguments to show that we can indeed iterate the procedure until all branches end in multigraphs that satisfy a halting criterion.

Finally, in Part 3, we verify that the algorithm with the properties given in Parts 1 & 2 does indeed yield the claimed upper bound. The idea is simple, even though the set-up in Parts 1 & 2 is messy. Our construction at each iteration generally yields (except in Criterion 2C, where we need a minor adjustment) multigraphs with cost at most 1, and we show that it guarantees a reduction of $\min\{n/3, \phi\}$ by a fixed factor, yielding the desired result $\chi'(G) \leq \phi(G) + \log\min\{n/3, \phi(G)\}$ through repeated application of Lemma A'.

### 3. Additional Background

Two classic 1-factor results will help assure that 1-factors or near 1-factors exist in multigraphs at key stages of our reduction.

**Lemma C.** (Tutte's 1-factor Theorem). Let G be a multigraph. There is a 1-factor of G if and only if for each set K of vertices of G, the number of odd-order components of G - K does not exceed |K|.

**Lemma D.** (Seymour [11]) Any r-graph has a 1-factor.

We will require also the following background Lemmas. They are referenced in [10], though appeared first in [6] and [11].

**Lemma E.** Let H be a nontrivial odd order subgraph of G. If H is $\Delta(G)$-overfull then $|\partial(V(H))| < \Delta(G)$; if H is $\Delta(G)$-full then $|\partial(V(H))| \leq \Delta(G)$.

**Lemma F.** Let G be r-regular. Then G is an r-graph if and only if $|\partial(S)| \geq r$ for each odd order subset S of V(G).

We now introduce some additional terminology. If k is any integer and H is a nontrivial odd-order multigraph, the *k-excess of H* (denoted ex(H, k)) is given by $e(H) - [k(n(H) - 1)/2]$; thus the k-excess of H gives the number of edges by which H is k-overfull. Note that ex(H, k) may be negative, and in fact ex(H, k) is 0 (positive) if and only if H is k-full (overfull). The *k-slack* of H (denoted sl(H, k)) is given by $(k + 1)[(n(H) - 1)/2] - e(H)$. Thus sl(H, k) gives the number of additional edges that H would need in order to (k + 1)-full ; note also that sl(H, k) + ex(H, k) = (n(H) – 1)/2. Often when we discuss k-excess and k-slack, we are considering H as a subgraph of a multigraph G, and use the value $k = \Delta(G)$. In that default case, we will sometimes use the simpler notation ex(H) and sl(H).

In the proof we will require a careful count of the effect that matching removal and splitting have on the excess of odd-order subgraphs. The following two technical observations will be of use throughout the proof.

**Observation 1.** (Effect of 1-factor removal on excess)

Let $G$ be a multigraph of even order $n$, and let $S$ be a nontrivial odd-cardinality subset of vertices of $G$. If $F$ is any 1-factor of $G$ and $\Delta = \Delta(G)$, then
$ex(<S>, \Delta-1; G-F) \leq ex(<S>, \Delta; G) + \min\{(n-|S|-1)/2, (|S|-1)/2\}$.

**Proof.** Since $\Delta = \Delta(G-F) + 1$, it follows that
$ex(<S>, \Delta-1; G-F) = ex(<S>, \Delta; G) + (|S|-1)/2 - r$, where $r$ gives the number of edges of $F$ which are incident with two vertices of $S$. But $F$ has $n/2$ edges, and at most $n - |S|$ of them do not have both incident vertices in $S$. Thus $r \geq |S| - n/2$, and the result follows. ◊

**Observation 2.** (Effect of shrinking on excess)

Let $G$ be a multigraph, and suppose that $S$ is a nontrivial odd-cardinality subset of $V(G)$. Let $R'$ be a nontrivial odd-cardinality set of vertices in $G_S$ and let $s$ be the vertex of $G_S$ that replaces $S$ upon shrinking $S$. For any integer $k$, if $s \notin R'$, then $ex(<R'>, k; G_S) = ex(<R'>, k; G)$, and
if $s \in R'$, then $ex(<R'>, k; G_S) = ex(<R>, k; G) - ex(<S>, k; G)$, where $R = R' - \{s\} \cup S$.

**Proof.** To find $ex(<R'>, k; G_S)$ we need only check to see how the numbers of edges and vertices in $<R'>$ are changed by the shrinking procedure. If $s \notin R'$, the result is clear, so assume $s \in R'$. Let $R$ be the set of vertices $R' - \{s\} \cup S$ in $G$. Then $|R| = |R'| + |S| - 1$, while
$|E(<R'>; G_S)| = |E(<R>; G)| - |E(<S>; G)|$.
Thus $ex(<R'>, k; G_S) = |E(<R'>; G_S)| - k(|R'|-1)/2$
$= |E(<R>; G)| - |E(<S>; G)| - k(|R'|-|S|+|S|-1)/2$
$= [|E(<R>; G)| - k(|R'|+|S|-2)/2] - |E(<S>; G)| + k(|S|-1)/2$
$= ex(<R>, k; G) - ex(<S>, k; G)$, as desired. ◊

In the following Lemma, note that because $G$ is $\Delta$-regular, $\Gamma(G) \geq \Delta$. It is possible, though, that $|S| = |V(G)| - 1$, in which case we would get $G_{S^c}$ to be isomorphic to $G$.

**Lemma G.** (Effects of shrinking maximum excess subgraph in a regular multigraph)

Let $G$ be a $\Delta$-regular multigraph of even order and suppose $<S>$ is an odd-order induced $\Delta$-full or $\Delta$-overfull subgraph which has maximum $\Delta$-excess in $G$ among all such subgraphs. Then
$\Delta(G_S) \leq \Delta$ and $\Gamma(G_S) \leq \Gamma(G)$, and similarly
$\Delta(G_{S^c}) \leq \Delta$ and $\Gamma(G_{S^c}) \leq \Gamma(G)$.

Proof. Because $\Delta(G) = \Delta$, and $<S>$ is $\Delta$-full or $\Delta$-overfull, $|\partial(S)| = |\partial(S^c)| \leq \Delta$. Thus we have both $\Delta(G_S) \leq \Delta$ and $\Delta(G_{S^c}) \leq \Delta$. Moreover, because $<S>$ and therefore also $<S^c>$ (because G is $\Delta$-regular) have maximum $\Delta$-excess, the other two inequalities in the statement of the Lemma follow, because by Observation 2 any subgraph containing the vertex from shrinking in $G_S$ or $G_{S^c}$ will have non-positive $\Delta$-excess. ◊

Our final preliminary result is similar in nature to Lemma G, but has the advantage that it can be applied in situations where $G$ is not regular.

**Lemma H.** (Effects of shrinking a minimum slack subgraph)

Let $G$ be a multigraph of even order, with $\Delta(G) < \Gamma(G) \leq \Delta(G) + 1$. Let $\Delta = \Delta(G)$, let $<S>$ be an odd order subgraph that is $\Delta$-overfull in $G$, with $sl(<S>, \Delta)$ a minimum among all such $\Delta$-overfull subgraphs. Then:

$\Delta(G_S) \leq \Delta(G)$, and $\Gamma(G_S) \leq \Delta(G) + 1$, and
$\Delta(G_{S^c}) \leq \Delta(G)$, and $\Gamma(G_{S^c}) \leq \Delta(G) + 1$.

**Proof.** As $<S>$ is k-full or k-overfull in $G$, by Lemma E we have $|\partial(S; G)| = |\partial(S^c; G)| < \Delta(G)$, so that $\Delta(G_S) \leq \Delta(G)$ and $\Delta(G_{S^c}) \leq \Delta(G)$, as desired.

Now let $M$ be a matching of new edges, its vertices identified with those of $S$, such that $|E(M)| = sl(<S>, \Delta)$ (so $t(<S> \cup E(M)) = \Delta(G) + 1$). Since adding a matching can increase $t(H)$ for any subset $H$ by at most 1, and since $S$ has minimum $\Delta$-slack among all overfull subgraphs of $G$, $\Gamma(G \cup E(M)) = \Delta(G) + 1$. By Lemma C, we can add edges to expand $G \cup E(M)$ to a multigraph $G^+$ which is a $(\Delta(G) + 1)$-graph. Then in $G^+$, $<S>$ is $(\Delta + 1)$-full and has maximum $(\Delta + 1)$-excess equal to 0, as does $<S^c>$. Thus by Lemma G, $\Gamma(G^+_S) \leq \Delta(G) + 1$ and $\Gamma(G^+_{S^c}) \leq \Delta(G) + 1$, so $\Gamma(G_S) \leq \Delta(G) + 1$ and $\Gamma(G_{S^c}) \leq \Delta(G) + 1$. ◊

## 4. Main Theorem

We are now ready to take on the proof of Theorem 1.

**Theorem 1.** For any multigraph G with order $n \geq 3$ and at least one edge,

$\chi'(G) \leq \phi(G) + \log_{3/2} (\min \{ n/3, \phi(G) \})$ if $n$ is even, and

$\chi'(G) \leq \phi(G) + \log_{3/2} (\min \{ (n+1)/3, \phi(G) \})$ if $n$ is odd.

Proof. As noted earlier, it suffices to prove the result for $n$ even, and by Lemma B we may also assume that G is an r-graph, where $r = \Delta(G) = \phi(G)$. Throughout the proof, we omit the log base, as it will always be 3/2. Because $n \geq 3$ and G has an edge, $\log (\min \{ n/3, \phi(G) \}) \geq 0$. The result is thus true when $n \leq 8$ for then $\chi'(G) = \phi(G)$ by Theorem C. The result also holds if $\phi(G) \leq 2$ because then each component of G is a path or even cycle and $\chi'(G) = \phi(G)$. Again the result is true if $\phi(G) = 3$ and $n > 8$, for in this case $\Delta(G) \leq 3$ so $\chi'(G) \leq 4$ by Shannon's Theorem, and $\log (\min \{ n/3, \phi(G) \}) \geq 1$. So assume that G has even order $n \geq 10$, that $\phi(G) \geq 4$, and the result is true for all multigraphs H with $n(H) \leq n(G)$ and $\phi(H) < \phi(G)$, or with $n(H) < n(G)$ and $\phi(H) \leq \phi(G)$.

*Part 1 – Algorithm Stage 1.* We first show that in Stage 1 we can construct a decomposition tree in accordance with the earlier discussion in which each terminal vertex multigraph H satisfies one of the Stage 1 halting criteria, without incurring cost more than 1.

Recall that our multigraph G has even order $n$ and is an r-graph, with $r = \Delta = \Delta(G) = \phi(G)$. By Lemma D, G has a 1-factor F. Clearly G - F is ($\Delta$-1)-regular; if $\phi(G - F) = \Delta - 1$ (so the cost of the reduction is 0), then

$\chi'(G) \leq 1 + \chi'(G - F) \leq 1 + \phi(G - F) + \log (\min \{ n/3, \phi(G - F) \})$
$\leq \phi(G) + \log (\min \{ n/3, \phi(G) \})$ so we are done.

So assume that G – F contains a ($\Delta$-1)-overfull subgraph.

Let $<S>$ be an induced ($\Delta$-1)-overfull odd order subgraph of $G - F$ which has maximum excess among all such subgraphs; as $G - F$ is regular, $ex(<S>, \Delta-1) = ex(<S^c>, \Delta-1)$. By Lemma G, $(G-F)_S$ and $(G-F)_{S^c}$ both have maximum degree $\Delta -1$, and we have both $\Gamma((G-F)_S) \leq \Delta$ and $\Gamma((G-F)_{S^c}) \leq \Delta$.

Let $H$ denote either $(G-F)_S$ or $(G-F)_{S^c}$. For the arguments below we will assume that $H = (G-F)_S$; the argument is essentially identical if we have $H = (G-F)_{S^c}$. Recall that $n$ denotes $n(G)$. Currently we have cost($G \to H$) = 1 because $\phi(H) = \phi(G)$ and mr($G \to H$) = 1. We note that one vertex $s$ in $H$ (the vertex replacing the shrunken set $S$) has degree less than $\Delta$-1 (by Lemma E), and all other vertices of $H$ have degree $\Delta$-1.

**Claim 1.** Any ($\Delta$-1)-overfull induced odd subgraph $<R>$ of $H$ has $ex(<R>, \Delta-1 ; H) \leq (n(G) - |R| - 1)/2$.

**Proof of Claim 1.**
If $s \notin R$, then $ex(<R>, \Delta-1; H) = ex(<R>, \Delta-1; G-F) \leq (n - |R| - 1)/2$, by Observation 1.
If $s \in R$, let $Q$ denote the subset of $V(G)$ given by $R - \{s\} \cup S$.
Then by Observation 2,
$ex(<R>, \Delta-1; H) = ex(<Q>, \Delta-1; G-F) - ex(<S>, \Delta-1; G-F)$
$\leq 0$ by the maximality of $ex(<S>, \Delta-1 ; G - F)$.
This yields the desired result. ◊

Let $\Delta^*$ denote the maximum degree of $H$. At this point, we have an even order multigraph $H$ with the following properties:
    (i) all vertices in $H$ have degree $\Delta^* = \Delta(G) -$ mr($G \to H$), except one vertex $s$ may have degree less than $\Delta^*$,
    (ii) $\phi(H) \leq \Delta^* + 1$, and so cost($G \to H$) $\leq 1$, and
    (iii) any induced $\Delta^*$-overfull subgraph $<R>$ of $H$ has $\Delta^*$-excess at most $(n(G) - |R| - 1)/2$.

The following Proposition shows that we can iterate our operations until all multigraphs obtained satisfy a Stage 1 Halting Criterion.

**Proposition 1.** Let H be a multigraph which does not satisfy a Stage 1 Halting Criterion, but does satisfy conditions (i), (ii), and (iii) above. Then either:

(a) H contains a 1-factor whose removal results in a new multigraph that satisfies properties (i) through (iii), or
(b) H contains a non-trivial splitting such that the resulting multigraphs either satisfy a Stage 1 Halting Criterion, or satisfy (i) though (iii), except that two vertices may have degree less than the others.

We will perform a sequence of 1-factor removals and splittings, using Proposition 1 to guarantee that we retain conditions (i) – (iii) and can therefore continue until all branches reach a Halting Criterion. If as in (b) we get two vertices with lower degree, we show in the argument following the proof of Claim 8 that we can add edges between the vertices to keep the desired properties (i) – (iii). Proposition 1 is proved below, but we first give the Halting Criteria and the basic algorithm for Stage 1 reduction, starting with multigraph H.

Any multigraph H obtained in Stage 1 will have cost(G→H) ≤ 1. Multigraph H obtained in our algorithm will be considered Stage 1 terminal (no further reduction needed in this stage) if it satisfies one of the following criteria.

**Stage 1 Halting Criteria for a multigraph H:**

Halting Criterion 1A. $\Delta(H) \leq 2$ or H has order at most 8.
   Note that by Theorem C, we know that then $\chi'(H) = \phi(H)$.

Halting Criterion 1B. Let s be the vertex with minimum degree in H. Stop if $ex(H-s, \Delta(H)) \geq n(H)/4$.
   We shall show that if we reach this halting criterion in our procedure, then two important relations hold:
   (a) $n(H) \leq 2\, n(G)/3$, and
   (b) $\Delta(H) \geq n(H)/2$, so $\phi(H) \geq n(H)/2$.

Halting Criterion 1C. H has a vertex of degree 0. We shall show also that in this case, $\Delta(H) \leq n(G)/3$. Clearly here H will not have a perfect matching. We will then need to continue the reduction procedure by going to Stage 2.

Stage 1 Reduction Algorithm   Let $\Delta^* = \Delta(H)$.

> If in H - s there is no $\Delta^*$-overfull odd order induced subgraph <R> for which $sl(<R>, \Delta^*; H\text{-}s) < sl(H\text{-}s, \Delta^*)$, then remove any 1-factor F; by Claims 5-7 below, such a 1-factor exists and H-F will satisfy properties (i)-(iii) with $\Delta^*$ replaced by $\Delta^* -1$.
> 
> Otherwise, such a multigraph R exists, so perform the splitting described in Claim 8 below, obtaining multigraphs $H_R$ and $H_{R^c}$. Multigraph $H_{R^c}$ will satisfy conditions (i)-(iii). Also $H_R$ will satisfy them, except there may now be two vertices s and $s_R$ with degree less $\Delta(H_R)$. If this occurs, add edges of the form $s s_R$ to $H_R$ until either :
> 
> (1) s or $s_R$ has degree $\Delta^*$, or
> (2) there is a subgraph <W> containing both s and $s_R$ for which $t(<W>) = \Delta^*+1$.
> 
> Call the resulting multigraph $H^+$.
> 
> If (1) occurs first, $H^+$ satisfies (i) – (iii) by the discussion following Claim 8.
> If (2) occurs first, perform another splitting, to obtain multigraphs $H^+_W$ and $H^+_{W^c}$. By the discussion following the proof of Claim 8, $H^+_W$ satisfies (i) – (iii), while $H^+_{W^c}$ satisfies Halting Criterion 1B.

As all multigraphs finally obtained now satisfy either a Halting Criterion or conditions (i) – (iii), we can iterate this procedure to eventually obtain all terminal multigraphs satisfying Halting Criteria for Stage 1, when Stage 1 is completed.

To verify Proposition 1, let $H$ be any multigraph which satisfies the properties (i) – (iii) above, but which does not satisfy any of the Stage 1 Halting Criteria. Let $p$ denote $n(H)$, and recall that $n$ denotes $n(G)$. We prove the Proposition by a series of claims.

**Claim 2.** The $\Delta^*$-excess of $H - s$ is $(\Delta^* - \deg(s; H))/2$.

**Proof of Claim 2.** Clearly $H - s$ has $[(p-1)\Delta^* - \deg(s; H)]/2$ edges. Thus $ex(H – s, \Delta^*) = [(p-1)\Delta^* - \deg(s; H)]/2 - (p-2)\Delta^*/2 = [\Delta^* - \deg(s)]/2$, as desired. ◊

Because H does not satisfy Halting Criterion 1C, we have $\deg(s; H) > 0$. We consider two cases.

Case A. In $H - s$ there is no $\Delta^*$-overfull odd order induced subgraph <R> for which $sl(<R>, \Delta^*; H-s) < sl(H-s, \Delta^*)$. We show that $H$ contains a 1-factor $F$ and $H – F$ satisfies (i) – (iii).

**Claim 3.** In Case A, if $<R>$ is $\Delta^*$-overfull in $H - s$, then $|R| \geq p/2$.
**Proof of Claim 3.** Assume $<R>$ is $\Delta^*$-overfull in $H - s$, and suppose that $|R| < p/2$, so $|R| \leq p/2 - 1$, since $p$ is even. Then $sl(<R>, \Delta^*) \leq p/4 – 1$ because $<R>$ is $\Delta^*$-overfull. But $sl(H – s, \Delta^*) = p/2 - 1 – ex(H – s, \Delta^*)$, and because $H$ does not satisfy Halting Criterion 1B, $ex(H-s, \Delta^*) < p/4$. Thus $sl(H – s, \Delta^*) > p/4 - 1 \geq sl(<R>, \Delta^*)$, a contradiction to our being in Case A. ◊

**Claim 4.** In Case A, for any nontrivial odd subset $R$ of $V(H)$ containing $s$, $<R; H>$ is not $\Delta^*$-overfull.
**Proof of Claim 4.** Suppose $R$ is an odd order subset of $V(H)$ containing $s$, and that $<R ; H>$ is $\Delta^*$-overfull. Then
$|\partial(R ; H)| < \Delta^* - (\Delta^* - \deg(s)) = \deg(s)$, so
$ex(<R^c>, \Delta^* ; H ) \geq (\Delta^* - \deg(s))/2 = ex(H – s, \Delta^*)$ by Claim 2.

But $R^c$ is a proper subset of H - s, so then
sl($<R^c>$, $\Delta^*$) < sl(H-s, $\Delta^*$), a contradiction because we are in Case A. ◊

**Claim 5.** In Case A there exists a 1-factor F of H.
**Proof of Claim 5.** Let K be any subset of V(H), and suppose the number of odd components of H - K is t > |K|. All odd components of H - K, except possibly the one containing s and at most one which is overfull (and therefore by Claim 3 contains at least half the vertices of H) have at least $\Delta^*$ edges in their coboundary in H. Since deg(s) > 0, by Claim 4 the odd component containing s would have at least 1 coboundary edge. Thus, the number of edges in the coboundary of K in H is at least (t-2)$\Delta^*$ + 1 ≥ |K|$\Delta^*$ + 1, since t and K have the same parity (because p is even). This yields a contradiction, as each vertex in K has degree at most $\Delta^*$. Thus, by Tutte's Theorem, H contains a 1-factor F. ◊

**Claim 6.** Let R be an odd order subset of V(H) such that $<R>$ is $\Delta^*$-overfull in H, and let F be a 1-factor of H. If we are in Case A, then
ex($<R>$, $\Delta^*$ - 1; H - F ) ≤ (n - |R| - 1)/2.
**Proof of Claim 6.** By Claim 4, R does not contain s. Also
ex($<R>$, $\Delta^*$ ; H) ≤ ex ($<H – s>$, $\Delta^*$ ; H) - (p - |R| - 1)/2, because we are in Case A. But by condition (iii), ex($<H – s>$, $\Delta^*$ ; H) ≤ (n-p)/2, so by Observation 1,
ex($<R>$, $\Delta^*$ - 1 ; H - F) ≤ ((n-p)/2 - (p -|R|-1)/2) + (p -|R|-1)/2 = (n – p)/2.
Thus ex($<R>$, $\Delta^*$ - 1 ; H - F) ≤ (n-p)/2 ≤ (n-|R| - 1)/2, as desired. ◊

**Claim 7.** For any 1-factor F of H, $\phi$(H - F) ≤ $\Delta^*$.
**Proof of Claim 7.** Clearly $\Delta$(H - F) = $\Delta^*$ - 1. We need to show that if R is any nontrivial odd order subset of V(H), then t($<R>$ ; H - F) ≤ $\Delta^*$ .
If t($<R>$ ; H) ≤ $\Delta^*$, this is certainly true. So, assume that $<R>$ is $\Delta^*$-overfull in H. Because H does not satisfy Halting Criterion 1B, ex(H – s, $\Delta^*$ ; H) < p/4, so ex( H – s, $\Delta^*$ ; H) ≤ (p – 2)/4 . Therefore, because we are in Case A,
ex($<R>$, $\Delta^*$ ; H) ≤ (p – 2)/4 – ( p – 1 - |R|)/2. Thus by Observation 1,
ex($<R>$, $\Delta^*$-1 ; H – F) ≤ (p - 2)/4 – (p – 1 - |R|)/2 + (p – 1 - |R|)/2 = (p – 2)/4.
But by Claim 3, |R| ≥ p/2, and thus t($<R>$ ; H - F) ≤ $\Delta^*$. ◊

Thus in Case A, by Claim 5 H contains a 1-factor F. Because H satisfies (i), all vertices in H - F have degree $\Delta(G) - mr(G \to (H - F))$, except possibly one vertex s of smaller degree. By Claim 7, $\phi(H - F) \leq \Delta(H - F) + 1$, so $cost(G \to H - F) \leq 1$. And by Claim 6, any induced overfull subgraph of H – F has $(\Delta^*-1)$-excess at most $(n - |R| - 1)/2$. Thus with the substitutions $H \leftarrow H - F$, $\Delta^* \leftarrow \Delta^* - 1$, the key properties (i) – (iii) remain satisfied.

Case B. There exists a $\Delta^*$-overfull odd order subgraph <R> of H - s such that $sl(<R>, \Delta^*) < sl(H - s, \Delta^*)$. We show that there is a splitting in which the resulting multigraphs (almost) satisfy (i) – (iii).

Let R be a subset of V(H) such that <R> has minimum $\Delta^*$-slack among all $\Delta^*$-overfull subgraphs of H. By Lemma H, $H_R$ and $H_{R^c}$ each have maximum degree at most $\Delta^*$, and $\phi$-values which are at most $\Delta^*+1$, so they satisfy property (ii). We wish to show that they also have property (iii).

**Claim 8.** In both $H_R$ and $H_{R^c}$, any nontrivial induced odd order subgraph <Q> has $\Delta^*$-excess at most $(n - |Q| - 1)/2$.

**Proof of Claim 8.**
First consider $H_R$ and let $s_R$ be the vertex of $H_R$ that was formed by shrinking R. Let <Q> be any nontrivial odd order induced subgraph of $H_R$.
If $s_R \notin Q$, then $ex(<Q>, \Delta^*; H_R) = ex(<Q>, \Delta^*; H) \leq (n - |Q| - 1)/2$, because condition (iii) is satisfied within H.
If $s_R \in Q$, then $ex(<Q>, \Delta^*; H_R) < ex(<Q - s_R \cup R>, \Delta^*; H)$ (by Observation 2, since <R> is $\Delta^*$-overfull)

$$\leq [n - (|Q| - 1 + |R|) - 1]/2 \text{ (by (iii) in H)}$$
$$< (n - |Q| - 1)/2, \text{ as desired.}$$

Now consider $H_{R^c}$, and let $s_{R^c}$ denote the vertex that replaces $<R^c>$ upon shrinking. Let <Q> be any nontrivial odd order induced subgraph of $H_{R^c}$.
If $s_{R^c} \notin Q$, then as before
$ex(<Q>, \Delta^*; H_{R^c}) = ex(<Q>, \Delta^*; H) \leq (n - |Q| - 1)/2$.
So, assume that $s_{R^c} \in Q$.

If $s \notin R^c$, then $<R^c>$ must be overfull in $H$, as $ex(<R^c>; H) \geq ex(<R>; H)$.
Thus $ex(<Q>, \Delta^*; H_{R^c}) \leq ex(<Q - s_{R^c} \cup R^c>, \Delta^*; H) \leq$

$[n - (|Q| - 1 + |R^c|) - 1]/2 < (n - |Q| - 1)/2$, as desired.

Finally, assume that $s_{R^c} \in Q$ and $s \in R^c$. Note that $s_{R^c}$ is the only possible vertex with degree less than $\Delta^*$ in $H_{R^c}$.
Suppose that $ex(<Q>; H_{R^c}) > (n - |Q| - 1)/2$. Then
$ex(<Q^c>; H_{R^c}) \geq ex(<Q>; H_{R^c}) > (n - |Q| - 1)/2 \geq (|Q^c| - 1)/2$ (the first inequality because $s_{R^c} \in Q$, the third because $|Q| + |Q^c| \leq n$). Thus
$t(<Q^c>; H_{R^c}) > \Delta^* + 1$, yielding a contradiction. Claim 8 now follows. ◊

Thus both $H_R$ and $H_{R^c}$ satisfy conditions (i), (ii), and (iii), with one possible exception: in $H_R$ there may be two vertices with degree less than $\Delta^*$: $s$ and $s_R$. If this occurs, we add parallel edges of the form $s s_R$ to $H_R$ until one of two things happen: either $s$ or $s_R$ has degree $\Delta^*$, or until there is an odd order subgraph $<W>$ containing both $s$ and $s_R$ for which $t(<W>) = \Delta^* + 1$. Let $H^+$ denote $H_R$ with the added edges. We first show that $H^+$ satisfies condition (iii). To do so, we show that as we were adding edges between $s$ and $s_R$ to $H_R$, we kept the property that any induced overfull odd subgraph $<Q>$ has $ex(<Q>, \Delta^*) \leq (n - |Q| - 1)/2$. If not, then since $H_R$ originally satisfied property (iii), both $s$ and $s_R$ must be in Q, so
$ex(<Q^c>; H^+) \geq ex(<Q>; H^+)$. But if $ex(<Q>) > (n - |Q| - 1)/2$, we have
$ex(<Q^c>) > (n - |Q| - 1)/2 \geq (|Q^c| - 1)/2$, so that $t(<Q^c>) > \Delta^* + 1$, a contradiction. Thus, Property (iii) is retained.

Thus in the first case above that completes $H^+$ (s or $s_R$ has degree $\Delta^*$ in $H^+$), $H^+$ now also satisfies conditions (i) – (iii). In the latter case (multigraph $W$ containing both $s$ and $s_R$ has $t(<W>) = \Delta^* + 1$), we use the shrinking operation to form multigraphs $H^+_W$ and $H^+_{W^c}$. Because
$t(<W>) = \Delta^* + 1$, $H^+_{W^c}$ satisfies Halting Criterion 1B, and so is terminal – it has two vertices with degree less than $\Delta^*$ but we have made enough progress. It is straightforward to check that it also contains the additional properties of Criterion 1B (because $ex(<W^c>) > ex(<W>)$, but $t(<W^c>) \leq t(<W>)$, we have

$|W| < |W^c|$, so $n(H^+_{W^c}) < 2n(G)/3$; also as $t(<W>) = \Delta^*+1$, $\Delta(H^+_{W^c}) \geq \Delta^*$).
It is also easy to check that the multigraph $H^+_W$ satisfies each of the key properties (i), (ii) and (iii). Thus the algorithm can continue.

We have thus verified Proposition 1, and shown that the reduction procedure can continue until the multigraphs obtained each satisfy one of the Halting Criterion 1A, 1B or 1C. To see that the Criteria 1B and 1C carry with them all the claimed properties, let $\Delta = \Delta(H)$.

First suppose that as in 1B we have $ex(H-s, \Delta) \geq n(H)/4$. A simple count shows that this forces $\Delta(H) \geq n(H)/2$, and so of course $\phi(H) \geq n(H)/2$. Moreover, since by condition (iii), $ex(H-s, \Delta) \leq (n - n(H))/2$, we have $n(H)/4 \leq (n - n(H))/2$ so $n(H) \leq 2n/3$, as desired.

If instead H satisfies 1C but not 1B, then because $H - s$ is $\Delta(H)$-regular, it has excess $\Delta(H)/2$, so $\Delta(H) < n(H)/2$. Thus if $n(H) \leq 2n/3$, we get the desired inequality $\Delta(H) \leq n/3$. On the other hand, if $n(H) > 2n/3$, then again by property (iii), $ex(H-s, k) \leq (n - n(H))/2$, so $\Delta(H)/2 \leq (n - n(H))/2$ and so again $\Delta(H) \leq n/3$.

If the original maximum degree of G was sufficiently small and a multigraph is terminated by Criterion 1C, we would not yet be guaranteed significant reduction to continue the algorithm and achieve our desired result. Thus we require Stage 2 of the reduction.

*Part 2 – Algorithm Stage 2*

We now show that we can continue with a Stage 2 reduction that mimics, with some adjustments, the procedure in the previous stage, using splitting and removal of matchings in order to reduce the multigraphs involved. Now however the matchings we remove will sometimes need to be one edge short of being a 1-factor (when one vertex has degree 0) in which case not all vertex degrees are reduced. Typically, then, we will be working with a multigraph H in which each vertex has degree either d or d+1 for some d ( these vertices we call Type 0 and Type 1 vertices, respectively), except possibly for a single vertex s for which deg(s) ≤ d. In such a situation, we say that H – s is *d-slack dominant* if sl(H – s, d) < sl( R, d) for any odd-order, non-trivial proper subgraph R of H – s which has coboundary cardinality at most d. A key observation is that then we must have ex(H-s, d) > ex(R, d) + (n(H - s) – n(R))/2 for each such R. Each time we remove a near 1-factor, we will create an additional Type 1 vertex, and generally this will be the only way such vertices are created. We seek to reduce the multigraphs so they satisfy Halting Criterion 2A or 2B below, mimicking what we did in Stage 1. However, the existence of the Type 1 vertices makes things more complicated, and thus we also allow as successful terminal multigraphs those satisfying Criterion 2C or 2D. Recall that G denotes the multigraph we began with in Stage 1.

**Stage 2 Halting Criteria** (Employed after reaching criterion 1C ).

Criterion 2A.  $n(H) \leq 8$ or $\Delta(H) \leq 2$ ; $\text{cost}(G \to H) \leq 1$.

Criterion 2B.  $n(H) \leq (2/3) n(G)$, $\phi(H) \geq n(H)/2$, $\text{cost}(G \to H) \leq 1$.

Criterion 2C.  $n(H) \leq \phi(G)$, $n(H) \leq n(G)/3$, $\text{cost}(G \to H) \leq 2$.

Criterion 2D. $\phi(H) \leq 2 n(G)/9$, $\phi(H) \leq \phi(G)/2$, $\text{cost}(G \to H) \leq 1$.

The flow of these Criteria is summarized in Figure 3 on the next page.

If the number of Type 1 vertices is more than d, this will show we have removed a significant number of matchings, and in fact we will show that the multigraph satisfies Criterion 2D (this is why Splitting Lemma 3 and the Matching Lemma that follow assume that the number of Type 1 vertices is at most d). The following Proposition therefore forms the basis for completing Stage 2. We will show in our algorithm that if t(H - s) > d+1, then H satisfies a Halting Criterion, so we assume in the Proposition that t(H - s) ≤ d+1.

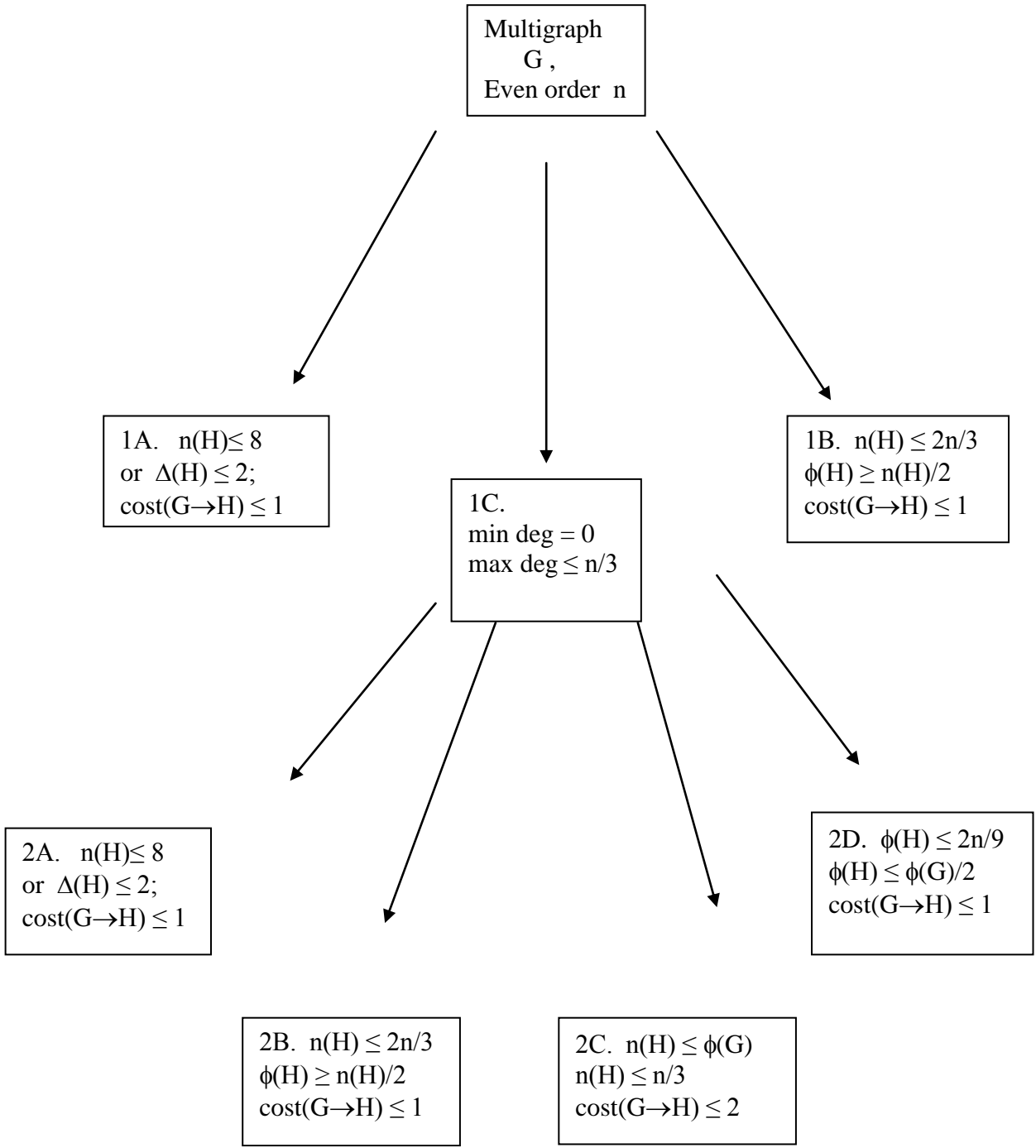

Figure 3. Flow of algorithm to halting criteria

**Proposition 2.**  Let H be an even order multigraph, with all vertex degrees either d or d+1, except one vertex s has degree (possibly 0) at most d.  Assume H has at most d vertices with degree d+1, that t(H - s) ≤ d+1 and ϕ(H) ≤ d+2. Then either

(a) there is a non-trivial splitting of H whose resulting multigraphs satisfy either some Stage 2 Halting Criterion or the same degree and ϕ-value restrictions as H (though there may be two vertices with degree less than d), or

(b) H has a perfect or (if deg(s) = 0) near-perfect matching whose removal leaves at most one additional vertex with maximum degree.

Using Proposition 2 as the foundation, we then show with some additional analysis that Stage 2 reduction can be completed as desired.

**Proof of Proposition 2.**  The proof will follow from a sequence of four lemmas; the first three provide splitting rules, the fourth guarantees an appropriate matching when needed.  Lemmas 1-3 handle situations where Γ(H) is respectively greater than, equal to, or less than d+1; these are used when H – s is not d-slack dominant.

**Splitting Lemma 1.**  Let H be a multigraph of even order, with Δ(H) ≤ d + 1, r vertices of degree d+1, and d+1 < Γ(H) ≤ d+2.  Let  < S >  be (d+1)-overfull in H, and assume that sl(< S >; d+1) is a minimum among all (d+1)-overfull subgraphs.  Then:

Δ($H_S$) ≤ d +1, Γ($H_S$) ≤ d+2, and $H_S$ has at most r vertices of degree d+1;

similarly,

Δ($H_{S^c}$) ≤ d +1, Γ($H_{S^c}$) ≤ d+2, and $H_{S^c}$ has at most r vertices of degree d+1.

**Proof.**  Because < S > is (d+1)-overfull, its coboundary cardinality is less than d+1, so shrinking S or $S^c$ will not create a new vertex of degree d+1, so the degree claims are easily verified.  Now let M be a matching of new edges, its vertices identified with those of S, such that | E(M) | = sl(< S >; d+1) (so t(< S > ∪ E(M)) = d+2).  Since adding a matching can increase t(R) for any subset R by at most 1, and since S has minimum (d+1)-slack among all (d+1)- overfull subgraphs of H, $\Gamma(H \cup E(M)) = d+2$. By Lemma B, we can add edges to expand $H \cup E(M)$ to a multigraph $H^+$ which is a (d+2)-graph. Then $<S>$ has maximum (d+2)-excess 0 in $H^+$, as does $<S^c>$. Thus by Lemma G, $\Gamma(H^+_S) \leq d + 2$ and $\Gamma(H^+_{S^c}) \leq d + 2$, so $\Gamma(H_S) \leq d+2$ and $\Gamma(H_{S^c}) \leq d+2$, and the result follows. ◊

**Splitting Lemma 2.** Let H be a multigraph of even order, with $\phi(H) \leq d + 1$, and suppose $<S>$ is (d+1)-full in H (i.e. $t(<S>) = d+1$). Let r denote the number of vertices of degree d+1 in H. Then:

$\Delta(H_S) \leq d +1$ and $\Gamma(H_S) \leq d+1$; moreover, $H_S$ has at most r vertices of degree d+1. Also, $\Delta(H_{S^c}) \leq d +1$, $\Gamma(H_{S^c}) \leq d+1$, and if $H_{S^c}$ has more than r vertices of degree d+1, then it has exactly r+1 vertices, each of degree d+1.

**Proof.** The bounds on the maximum degree and $\Gamma$ values follow readily from Lemma G in the usual fashion if we first add edges to make the multigraph a (d+1)-graph. The graphs obtained by the shrinking can have an additional vertex of degree d+1 only if the coboundary cardinality $|\delta(S)| = d+1$ in H. But then in order to have $t(<S; H>) = d+1$, all vertices in S must have degree d+1. The result now follows. ◊

**Splitting Lemma 3.** Let H be a multigraph with even order $p \geq 4$. Suppose $d \geq 2$, and that except for one vertex s for which $\deg(s) \leq d$, each vertex has degree d or d+1, with $r \leq d$ vertices having degree d+1. Suppose that H – s is not d-slack dominant; that is, there exists a non-trivial odd order proper subgraph of H – s with coboundary cardinality at most d and whose d-slack is less than or equal to $sl(H - s, d)$. If $\Gamma(H) \leq d + 1$ then there is an odd-order subset S of V(H) with $1 < |S| < p-1$ such that

$\Delta(H_S) \leq d+1$, $\Gamma(H_S) \leq d + 1$, and the number of vertices with degree d+1 in $H_S$ is at most r, and similarly

$\Delta(H_{S^c}) \leq d+1$, $\Gamma(H_{S^c}) \leq d + 1$, and the number of vertices with degree $d+1$ in $H_{S^c}$ is at most $r$.

**Proof.**

Let R be a non-trivial odd order proper induced subgraph of H – s with coboundary cardinality at most d and sl( R, d) ≤ sl(H-s, d), with R chosen to have minimum slack among such subgraphs. Let M be a matching of new edges, its vertices identified with vertices of R which have degree d in H, such that | E(M) | = sl( R, d) (so t( R ∪ E(M)) = d + 1). Such a matching is guaranteed to exist because R is d-full or d-overfull with | cobound(R) | ≤ d and each vertex of R has degree at least d in H. Since adding a matching can increase t(L) for any subgraph L by at most 1, and since R has minimum d-slack among subgraphs of H - s with coboundary cardinality at most d, Γ(H - s ∪ E(M)) = d + 1. It may be however that Γ(H ∪ E(M)) > d+1, if there were a subgraph of H ∪ E(M) containing s that is (d+1)-overfull. Thus, we select M′ to be a minimal subset of M such that H ∪ E(M′) has a (d+1)-full subgraph (it is quite possible that M′ = M), and let S be the vertex set of that subgraph. By construction, the subgraph of H induced by S has coboundary cardinality at most d (if M′ ≠ M and thus S contains s, its coboundary cardinality is less than d+1 because s has degree at most d, and < S; H ∪ E(M′)> is (d+1)-full). By Lemma B, we can add edges to expand H ∪ E(M′) to a multigraph $H^+$ which is a (d + 1)-graph. Then < S > has maximum excess 0 in $H^+$, as does < $S^c$ >. Thus by Lemma G, $\Gamma(H^+_S) \leq d + 1$ and $\Gamma(H^+_{S^c}) \leq d + 1$, so $\Gamma(H_S) \leq d + 1$ and $\Gamma(H_{S^c}) \leq d + 1$. Finally, as the subgraph of H induced by S has coboundary cardinality at most d, the degree conditions on $H_S$ and $H_{S^c}$ follow as desired. ◊

**Matching Lemma.** Let H be a multigraph with even order p ≥ 4. Suppose d ≥ 2, and that except for one vertex s for which deg(s) ≤ d, each vertex has degree d or d+1, with at most d vertices having degree d+1, and at least one vertex other than s having degree d. Suppose that ex(H-s, d) < p/4, and that H - s is d-slack-dominant. Then:
- if deg(s) > 0, then H has a perfect matching, and
- if deg(s) = 0, then there is a vertex v of H with degree d, such that H – s – v has a perfect matching.

**Proof.** Let H be a multigraph with the degree properties given above, and suppose that H – s is d-slack-dominant with ex (H-s, d) < p/4. Because the excess is an integer and p is even, then ex (H-s, d) ≤ (p-1)/4 .

Let H' be that multigraph obtained from H according to the rule that H' = H if deg(s) >0, and if deg(s)=0 then H' = H + e, where e is an edge added between s and some vertex v with degree d in H. If some vertex of degree d is adjacent to a vertex of degree d+1 in H, we choose v to be such a vertex; otherwise we place no additional restriction on the choice of v.
Note that H' – s = H - s, so ex(H' - s, d) ≤ (p-1)/4 and H'-s is d-slack-dominant. However, H' has up to d+1 vertices of degree d+1. The Lemma will be proved if we can show that H' has a perfect matching.

Suppose that H' has no perfect matching. Because H' has even order, by Tutte's Theorem there is a set S of m vertices in H' such that H' – S has at least m+2 odd order components. We consider two cases.

Case 1. A non-trivial odd component R of H'–S not containing s has coboundary cardinality k ≤ d in H'.

Because ex(H'-s, d) ≤ (p-1)/4, and H' - s is d-slack-dominant, at most one such component R can exist, for it must have order greater than p/2. Let x , y respectively denote the number of vertices of H'- s not in R which have degree d, d+1 in H', and let j = deg(s). Then
ex(R, d) + (y + k – j)/2 = ex(H'- s, d) > ex(R, d) + (x + y)/2, the inequality following as H'-s has x+y more vertices than does R and sl(R, d) > sl(H' – s, d) because H' - s is d-slack dominant. Thus k – j > x.

Clearly the m vertices of S have degree sum at most m(d+1) in H'. But now consider the m+2 or more odd order components of H' – S. In addition to the component containing s (if there is one), and the component R, let a of them have d edges in H' joining them to members of S (each of these must be singleton vertices by the uniqueness of R), and the remaining b components have at least d+1 edges each joining them to members of S, where a+b ≥ m.

But $a \le x < k-j$. Thus the total number of edges in $H'$ joining the odd components of $H' - S$ to vertices in $S$ is at least
$k + ad + b(d+1) \ge k + b + md > a + b + md \ge m(d+1)$, a contradiction since the m vertices of S each have degree at most d+1. Thus $H'$ must have a perfect matching.

Case 2. Each non-trivial odd component of $H' - S$ not containing s has coboundary cardinality at least d+1 in $H'$.

Because the degree of s is at least 1, the (possibly trivial) odd components of $H' - S$ have total coboundary cardinality at least $(m+1)d + 1$ in $H'$. On the other hand, letting t denote the number of Type 1 vertices (vertices with degree d+1) in S, the degree sum in $H'$ of vertices in S is at most $md + t \le md + d + 1 = (m+1)d + 1$. Because this sum must be at least as large as the coboundary cardinality of the odd components of $H' - S$, we conclude we must have equality in the equation above so that $t = d + 1$, so S contains all d+1 of the Type 1 vertices of $H'$ (one of which is the earlier designated vertex v incident with e), with all their incident edges joining them to odd component vertices of $H' - S$. But v is then not adjacent to a vertex of degree d+1, contradicting the rule for its selection. The result follows. ◊

The four previous Lemmas together provide the proof of Proposition 2.

Finally, we present our reduction algorithm for Stage 2. Denote by $G_2$ the multigraph with which we enter Stage 2. Recall that $\Delta(G_2) \le n/3$, where $n = n(G)$. After a sequence of reductions in which m matchings have been removed, we let $d = \Delta(G_2) - m$. Thus d would be the maximum degree of our current multigraph H, if each matching had been a 1-factor. Suppose our multigraph H has $\Delta(H) \le d+1$ and $\phi(H) \le d+2$, that vertex s has minimum degree, and H does not satisfy any of the Criteria 2A, 2B, 2C or 2D.

The next step of the algorithm is given; we would iterate on any resulting multigraphs until all reach one of the four Stage 2 Halting Criteria.

1. If $\phi(H) > d+1$, let $<S>$ have minimum $(d+1)$-slack value among $(d+1)$-overfull subgraphs, and perform the splitting to get reduced subgraphs $G_S$ and $G_{S^c}$ as in Splitting Lemma 1. (By observation D below, this splitting in non-trivial and results in two multigraphs whose order is less than $n(H)$).

2. Else if $\phi(H) = d+1$, and there exists a subgraph $<S>$ with $1 < |S| < n(H) - 1$ and $t(<S>) = d+1$, perform the splitting described in Splitting Lemma 2. We note that if the number of Type 1 vertices increases in one of these shrinkings, the resulting multigraph has all vertices of degree $d+1$, and it is easy to check that this is then a terminal multigraph satisfying Criterion 2B, as shown in observation B below.

3. Else if there are at least two vertices $s, w$ with degree less than $d$, add an edge $sw$;
   Note: if we reach this operation, adding the edge will not increase $\phi$ above $d+1$, else we would be in Case 2.

4. Else $H$ satisfies the conditions of either Splitting Lemma 3 or the Matching Lemma, depending on whether on not $H - s$ is $d$-slack dominant. We apply the appropriate Lemma to get the desired reduction.

At each step we reduce either $d(H)$ or $n(H)$ (or add an edge without increasing them), so the algorithm eventually terminates. We complete the verification of Proposition 2 and the Stage 2 algorithm by verifying the following Observations on a non-terminal multigraph $H$ formed during Phase 2.

A. $\phi(H) \leq d(H) + 2$. (This guarantees that $cost(G \rightarrow H) \leq 2$).
   To see this, note that by the algorithm design, the first time we obtain a multigraph $H$ for which $\phi(H) > d(H) + 1$, we have just removed a matching, and $\phi(H) = d(H) + 2$. But for a multigraph $H$ with these properties, Case 1 of the algorithm applies, and thus the $\phi$-value will never surpass the d-value by more than 2.

B. If a new Type 1 vertex is created in Splitting Lemma 2, the resulting multigraph $H_{S^c}$ is terminal. (Thus if the algorithm needs to continue, the only way that Type 1 vertices have been created is by removal of near 1-factors).

To see this, recall that in the algorithm there are only two ways that a Type 1 vertex can be created: using the Matching Lemma to remove a near 1-factor may change a Type 0 vertex to a Type 1 vertex, or the new vertex created in forming $H_{S^c}$ in Splitting Lemma 2 may be Type 1. So, assume that we reach the situation where a Type 1 vertex is produced by Splitting Lemma 2. By that Lemma, all vertices in $H_{S^c}$ have degree $d(H) + 1$. Thus

$n(H_{S^c}) \leq 1 + mr(G_2 \to H) \leq \Delta(G_2)$. Thus $n(H_{S^c}) \leq \phi(G)$, and $n(H_{S^c}) \leq n(G)/3$, so $H_{S^c}$ is terminal under criterion 2C.

C. If $\phi(H) \leq d(H)+1$, at most $d = d(H)$ vertices are Type 1 (this is needed to apply Splitting Lemma 3 or the Matching Lemma).

To verify this, recall that when an iteration on a multigraph that is not terminal is being considered, the number of Type 1 vertices is at most the matching reduction number $mr(G_2 \to H)$, measured starting from $G_2$. Suppose there are $r > d$ Type 1 vertices. As $G_2$ is the multigraph that began Stage 2, we have that $n/3 \geq \Delta(G_2) \geq d+r \geq 2d + 1$, and therefore $\phi(G_2) \geq 2d+2$ since $\phi(G_2) > \Delta(G_2)$. We then have $\phi(H) \leq \phi(G)/2$. Also, since $(n/3)/(2d+1) \geq 1$, $\phi(H) \leq [(d+1)/(2d+1)](n/3) < 2n/9$ for $d \geq 2$. Thus $H$ would satisfy Criterion 2D, a contradiction.

D. $t(H-s) \leq d+1$. (Thus if we do the splitting in Case 1 of the algorithm, we get a real reduction, as we are not shrinking on $H-s$ and its complement).

To verify this property, suppose $t(H-s) > d+1$, and again let $r$ denote the number of Type 1 vertices in H. As before, we must have $n/3 \geq \Delta(G_2) \geq d+r$, and also $\phi(G) \geq \phi(G_2) \geq d+r$. But $t(H-s) > d+1$ implies that if H has k vertices other than s with degree less than $d+1$, $H-s$ must have at least $(d+1)(r + k - 1)/2 + 1$ edges, so its degree sum, which can be no more than $r(d+1) + kd$, must be at least $(d+1)(r + k - 1) + 2$, so that $k \leq d-1$.
Thus $n(H) \leq r + d$. It follows that $n(H) \leq \phi(G)$, and $n(H) \leq n/3$, so that H would satisfy Halting Criterion 2C, a contradiction.

*Part 3 – Reduction Verification*

It remains to show that if each terminal multigraph H in an iteration of the reduction process always satisfies one of the Halting Criteria 1A, 1B, 2A, 2B, 2C or 2D, then the theorem is true. In view of Lemma A′, it suffices to show that when cost(G→H) = k, we have
min { n(H)/3, $\phi$(H) } ≤ $(2/3)^k$ min { n(G)/3, $\phi$(G) }.

We consider the possible criteria in turn. If H satisfies 1A or 2A, it is actually possible that min{ n(H)/3, $\phi$(H) } is not sufficiently reduced. But we know that $\chi'$(H) = $\phi$(H) under these conditions, so we can remove additional matchings that each time reduce $\phi$(H) by one, so that min{ n(H)/3, $\phi$(H) } is eventually small enough.

Next suppose H satisfies 1B or 2B. If min { n(G)/3, $\phi$(G) } = n(G)/3, the reduction in n(H) is sufficient, so assume min { n(G)/3, $\phi$(G) } = $\phi$(G). But n(H) ≤ 2 $\phi$(H) ≤ 2 $\phi$(G), so min { n(H)/3, $\phi$(H) } ≤ (2/3) min { n(G)/3, $\phi$(G) } as desired.

Next suppose H satisfies Criterion 2C. Again, if min { n(G)/3, $\phi$(G) } = n(G)/3, the reduction in n(H) is sufficient, so assume min { n(G)/3, $\phi$(G) } = $\phi$(G). Because n(H) ≤ $\phi$(G), we have
min { n(H)/3, $\phi$(H) } } ≤ $\phi$(G)/3 ≤ $(2/3)^2$ min { n(G)/3, $\phi$(G) }, as desired.

Finally, suppose H satisfies Criterion 2D. If min { n(G)/3, $\phi$(G) } = n(G)/3, then since $\phi$(H) ≤ 2n(G)/9, we have the desired relation
min { n(H)/3, $\phi$(H) } ≤ (2/3) min { n(G)/3, $\phi$(G) }. Finally, if
min { n(G)/3, $\phi$(G) } = $\phi$(G), then because $\phi$(H) ≤ $\phi$(G)/2, we again have the desired reduction.

Thus since we can achieve the reductions with Halting Criteria as claimed, the proof follows.

We close with three comments. First, it is interesting to note that interrelations between n and $\phi$ were a key to this proof, so that for example trying to directly prove that $\chi'(G) \leq \phi(G) + \log(\phi(G))$ does not appear to be any simpler than getting our result. Second, in order to improve the result in this paper significantly, the most promising approach would seem to be to be more selective in choosing matchings that are removed; in our process, the matchings were selected arbitrarily. Finally, in trying to implement the algorithm, the main stumbling point would be identifying subgraphs whose slack is a minimum; indeed, this task could be NP-hard. It seems natural to instead adjust the algorithm to shrink on subgraphs with maximum fractional chromatic index. Although this approach fits better with other theory, it creates other complexities that would need to be overcome.